\documentclass[a4paper]{amsart}

\usepackage{amssymb}

\usepackage{amscd}

\usepackage{amsthm}
\theoremstyle{plain}
\newtheorem{thm}[]{Theorem}
\newtheorem{lem}{Lemma}
\newtheorem{cor}[lem]{Corollary}
\theoremstyle{remark}
\newtheorem{rmk}[]{Remark}

\usepackage{hyperref}
\usepackage{mathscinet}
\usepackage[initials,nobysame,lite,msc-links]{amsrefs}

\newcommand{\ensemble}[1]{\left\{ #1 \right\}}
\newcommand{\suchthat}{\mid}
\newcommand{\norm}[1]{\left\| #1 \right\|}
\newcommand{\absolute}[1]{\left| #1 \right|}
\newcommand{\N}{\mathbb{N}}
\newcommand{\Z}{\mathbb{Z}}
\newcommand{\C}{\mathbb{C}}

\DeclareMathOperator{\Tr}{Tr}
\DeclareMathOperator{\tr}{tr}
\DeclareMathOperator{\Ad}{Ad}
\DeclareMathOperator{\Id}{Id}
\DeclareMathOperator{\sign}{sign}
\DeclareMathOperator{\PU}{PU}

\begin{document}

\title[symmetric group subfactors]{On subfactors arising from asymptotic representations of symmetric groups}
\author{Makoto Yamashita}
% \date{}  % Activate to display a given date or no date
\address{Graduate School of Mathematical Sciences, The University of Tokyo, 3-8-1 Komaba, Tokyo, 153-8914, JAPAN}
\email{makotoy@ms.u-tokyo.ac.jp}

\begin{abstract}
We consider the infinite symmetric group and its infinite index subgroup given as the stabilizer subgroup of one element under the natural action on a countable set.  This inclusion of discrete groups induces a hyperfinite subfactor for each finite factorial representation of the larger group.  We compute subfactor invariants of this construction in terms of the Thoma parameter.
\end{abstract}

\maketitle

\section{Introduction}

The unitary representations of the infinite symmetric group $S_\infty$ have exhibited an interesting interplay of combinatorics and various aspects of analysis, such as probability theory.  Thoma~\cite{MR0173169} gave a parametrization of the II$_1$ factorial representations of $S_\infty$ in terms of a sequence of real parameters.  Subsequently there appeared several alternative proofs discovered by Vershik--Kerov~\citelist{\cite{MR614033}\cite{MR639197}}, Okounkov~\cite{MR1283250}, and Hirai~\cite{MR2118037}.  There have been various developments on this subject by Okounkov~\cite{MR1691646}, Kerov--Okounkov--Olshanski~\cite{MR1609628}, Bo{\.z}ejko--Gu{\c{t}}{\u{a}}~\cite{MR1923173}, and Kerov--Olshanski--Vershik~\cite{MR2104794}, among others.

The representations of $S_\infty$ has also been an important example in the theory of operator algebras, as evidenced by the construction of a hyperfinite II$_1$ factor using its regular representation by Murray--von Neumann, which was among the pioneering works of the subject.  When one considers the structure of the representations of $S_\infty$ through operator algebras, the first important phenomenon is the uniqueness of the hyperfinite (AFD) II$_1$ factor.  It implies that all the II$_1$ factorial representations of the locally finite group $S_\infty$ generate isomorphic von Neumann algebras.

Regarding the structure of the hyperfinite II$_1$ factor, the theory of subfactors initiated by Jones~\cite{MR696688} was very fruitful and revealed many surprising combinatorial aspects of the structure of inclusion between such factors.  One should note that the infinite symmetric group also has a similar structure of inclusions, as it contains a subgroup isomorphic to itself, given as the stabilizer of a single point under the natural representation on a countable set. This inclusion of discrete groups induces an inclusion of II$_1$ factors for each II$_1$ factorial representation of $S_\infty$.

In this paper we investigate the hyperfinite subfactors obtained this way and investigate how subfactor invariants reflect the Thoma parameter.  First we obtain the characterization of irreducibility in terms of the Thoma parameter (Theorem~\ref{thm:irred-characterization}), next the characterization of being an infinite index inclusion as the faithfulness of the corresponding state on the group algebra (Theorem~\ref{thm:subf-infin-ind}), and then an estimate of the subfactor entropy (Theorem~\ref{thm:subf-entropy-estim}).

Finally, we note that Gohm--K\"{o}stler~\cite{arXiv:1005.5726} recently obtained the same characterization of irreducibility as well as the new proof of Thoma's classification result.

\section{Preliminaries}

For each positive integer $n$, we identify the $n$-th symmetric group $S_n$ with the group of the bijections on the set $\ensemble{0, 1, \ldots, n - 1}$.  These groups form an increasing sequence by letting the elements of $S_n$ act on $\ensemble{n, \ldots, m - 1}$ trivially when $n < m$.  Given a sequence $m_1, \ldots, m_n$ of distinct numbers, we write $(m_1, m_2, \ldots, m_n)$ for the cyclic permutation $s$ defined by $s(m_k) = m_{k+1}$ for $k = 1, \ldots, n - 1$ and $s(m_n) = m_1$.  
  
The infinite symmetric group $S_\infty$ is defined to be the union of the groups $S_n$ for $n \in \N_{>0}$.  It is identified with the set of the permutations of $\N = \ensemble{0, 1, \ldots}$ which move only a finite number of elements. It is a countable infinite conjugacy class group, meaning that the conjugacy classes are infinite sets except for the trivial class of the neutral element.

For each $n \in \N$, let $S_{n\le}$ denote the subgroup of $S_\infty$ consisting of the elements which fix the numbers in $\ensemble{0, \ldots, n - 1}$.  Thus we have a decreasing sequence of infinite groups $S_\infty \supset S_{1\le} \supset S_{2\le} \cdots$, all isomorphic to $S_\infty$ and having infinite index at each inclusion.

The unitary representations of $S_\infty$ whose image generate (hyperfinite) II$_1$ factors are classified by Thoma (in the following we follow the treatment of Wassermann \cite{wassermann-thesis}*{Chapter III}).  We briefly recall their classification.  They are parametrized by the Thoma parameter
\[
\kappa = (\alpha_{1}, \alpha_{2}, \ldots; \beta_1, \beta_2, \ldots; \gamma),
\]
 consisting of two nonincreasing sequences of nonnegative real numbers $\alpha_{1} \ge \alpha_{2} \ge \cdots$, $\beta_1 \ge \beta_2 \ge \cdots$, and another nonnegative number $\gamma$ satisfying 
\[
\sum_{i=1}^\infty \alpha_i + \sum_{i=1}^\infty \beta_i + \gamma = 1.
\]
Let $\kappa$ be a Thoma parameter as above.  Then we construct a positive definite function $\tau_\kappa$ on $S_\infty$ by putting
\[
\tau_\kappa(s) = \sum_{i=1}^\infty \alpha_i^n - (-\beta_i)^n
\]
when $s$ is a cyclic permutation of $n$ numbers, and setting $\tau_\kappa(s_0s_1\cdots s_m) = \prod \tau_\kappa(s_k)$ when $s_0, s_1, \ldots, s_m$ are cyclic permutations with nonintersecting supports.  Thus, under the above notation we have
\[
\tau_\kappa((m_1, \ldots, m_n)) =  \sum_{i=0}^\infty \alpha_i^n - (-\beta_i)^n
\]
for $n$ distinct integers $m_1, \ldots, m_n$.

The function $\tau_\kappa$ defines an extremal tracial state on the group C$^*$-algebra $C^*S_\infty$ of $S_\infty$, hence the associated Gelfand--Naimark--Segal representation of $C^*S_\infty$ generates a hyperfinite II$_1$ factor unless $\alpha_1 = 1$ or $\beta_1 = 1$.  Moreover, these representations of $S_\infty$ for different values of $\kappa$ are never mutually quasi-equivalent.

Let $M_\kappa$ denote the II$_1$ factor generated by $S_\infty$ under the GNS representation with respect to $\tau_\kappa$ and let $N_\kappa$ denote the subfactor of $M_\kappa$ generated by the image of $S_{1\le}$.  In the rest of the paper we analyze the inclusion of von Neumann algebras $N_\kappa \subset M_\kappa$.

\section{Characterization of Irreducibility}
 
 Let $T_n$ be the permutation group of the set $\ensemble{1, \ldots, n-1}$.  The group $S_{1\le}$ is the union of the increasing sequence of the groups $T_n$.  The conditional expectation onto the relative commutant of $T_n$ inside $M_\kappa$ can be computed using the following lemma.
 
\begin{lem}\label{lem:condi-exp-transp}
 Let $s$ be the image of the transposition $(0, 1)$ under the representation $\pi_\kappa$.  Then the image of $s$ under the conditional expectation $E_{\pi_\kappa(T_n)' \cap M_\kappa}^{M_\kappa}$ is equal to $ \frac{1}{n-1} \sum_{j=1}^{n-1}\pi_\kappa((0, j))$.
 \end{lem}
 
\begin{proof}
For an arbitrary element $x$ of $M_\kappa$, one has
 \[
 E_{\pi_\kappa(T_n)' \cap M_\kappa}(x) = \frac{1}{(n-1)!} \sum_{t \in T_n} \pi_\kappa(t) x \pi_\kappa(t^{-1}).
 \]
Under the conjugation action by $T_n$, the stabilizer subgroup of $(0, 1)$ is the subgroup of $T_n$ consisting of the permutations of the set $\ensemble{2, \ldots, n - 1}$, and the conjugates of $(0, 1)$ are $(0, j)$ for $j = 1, \ldots, n - 1$. Hence we obtain
 \[
 E_{\pi_\kappa(T_n)' \cap M_\kappa}(\pi_\kappa((0, 1))) = \frac{1}{n-1} \sum_{j=1}^{n-1}\pi_\kappa((0, j)),
 \]
 which proves the assertion.
 \end{proof}
 
Recall that a subfactor $N \subset M$ is said to be irreducible when the relative commutant $N' \cap M$ equals $\C$.  We have the following characterization of irreducibility for the inclusion $N_\kappa \subset M_\kappa$.

\begin{thm}\label{thm:irred-characterization}
 The following conditions are equivalent:
\renewcommand{\theenumi}{\roman{enumi}}
 \begin{enumerate}
\item \label{irred-phenom}The inclusion $N_\kappa \subset M_\kappa$ is irreducible.
\item \label{irred-par-condi}One of the following conditions on $\kappa$ holds:
\begin{enumerate}
\item\label{par-tensor-trace-case} \label{eq:thoma-par-itpfi-perm}
There exists an integer $n > 0$ such that
\begin{equation*}
\kappa = (1/n, \ldots, 1/n, 0, 0, \ldots; 0, 0, \ldots; 0).
\end{equation*}
\item\label{eq:thoma-par-itpfi-perm-det}
There exists an integer $n > 0$ such that
\begin{equation*}
\kappa = (0, 0, \ldots; 1/n, \ldots, 1/n, 0, 0, \ldots; 0).
\end{equation*}
\item\label{par-reg-rep-case} One has $\alpha_i = \beta_i = 0$ for any $i$ (hence $\gamma = 1$).
\end{enumerate}
\end{enumerate}
\renewcommand{\theenumi}{\arabic}
\end{thm}

\begin{proof}
 \eqref{irred-phenom} $\Rightarrow$ \eqref{irred-par-condi}: Suppose that the inclusion $N_\kappa \subset M_\kappa$ is irreducible.  As in Lemma~\ref{lem:condi-exp-transp}, let $s$ be the image of the transposition $(0, 1)$ in $M_\kappa$ and $E_{N_\kappa' \cap M_\kappa}$ be the conditional expectation of $M_\kappa$ onto $N_\kappa' \cap M_\kappa$. By the irreducibility assumption and the general principle $\tau_\kappa(E_{N_\kappa' \cap M_\kappa}(x)) = \tau_\kappa(x)$, the image of $s$ under $E_{N_\kappa' \cap M_\kappa}$ must be equal to the scalar $\tau_\kappa(s)$.
 
By Lemma~\ref{lem:condi-exp-transp}, we have
\[
E_{\pi_\kappa(T_n)' \cap M_\kappa}(s)^2 = \frac{ n - 1 + \sum_{1 \le i \neq j \le n - 1} \pi_\kappa((0, i, j))}{(n - 1)^2}.
\]
Hence we obtain
 \begin{equation}\label{eq:sq-norm-tn-rel-comm-exp-s}
 \norm{E_{\pi_\kappa(T_n)' \cap M_\kappa}(s)}_2^2 = \tau_\kappa(  E_{\pi_\kappa(T_n)' \cap M_\kappa}(s) ^2) =\frac{1 + (n-2)\tau_\kappa((0, 1, 2))}{n-1}.
 \end{equation}
 
 On the other hand, we have $E_{N_\kappa' \cap M_\kappa}(s) = \lim_{n\rightarrow\infty} E_{\pi_\kappa(T_n)' \cap M_\kappa}(x)$ in the $2$-norm topology.  Since $ \norm{E_{\pi_\kappa(T_n)' \cap M_\kappa}(s)}_2^2 $ converges to $\tau_\kappa((0, 1, 2))$ as $n \rightarrow 0$ by \eqref{eq:sq-norm-tn-rel-comm-exp-s}, one has to have $ \tau_\kappa(s)^2 = \tau_\kappa((0, 1, 2))$ that is equivalent to
 \begin{equation}\label{eq:s-tau-square-eq-3-elem-tau}
 \sum_{i=1}^\infty \alpha_i^3 + \sum_{i=1}^\infty \beta_i^3 = (\sum_{i=1}^\infty \alpha_i^2 - \sum_{i=1}^\infty \beta_i^2)^2.
 \end{equation}

Let $(\gamma_i)_{i \in \Z \setminus \ensemble{0}}$ be the sequence defined by $\gamma_i = \alpha_i$ for $i > 0$ and by $\gamma_i = \beta_{-i}$ for $i < 0$.  Then we have
 \begin{multline}\label{eq:3rd-mom-gerater-sq-2nd-mom}
 \sum_{i=1}^\infty \alpha_i^3 + \sum_{i=1}^\infty \beta_i^3 - (\sum_{i=1}^\infty \alpha_i^2 + \sum_{i=1}^\infty \beta_i^2)^2 \\
  = (\sum_{i=1}^\infty \gamma_i^3) (\sum_{i=1}^\infty \gamma_i) - (\sum_{i=1}^\infty \gamma_i^2)^2 = \sum_{i < j} \gamma_i \gamma_j (\gamma_i - \gamma_j)^2 \ge 0.
\end{multline}
Combining this with the inequality
\begin{equation}\label{eq:diff-sq-smaller-sum-sq}
 (\sum_{i=1}^\infty \alpha_i^2 + \sum_{i=1}^\infty \beta_i^2)^2 \ge (\sum_{i=1}^\infty \alpha_i^2 - \sum_{i=1}^\infty \beta_i^2)^2,
\end{equation}
it follows that \eqref{eq:s-tau-square-eq-3-elem-tau} holds if and only if the the above inequalities are actually equalities.

On one hand, in order to have the equality in \eqref{eq:3rd-mom-gerater-sq-2nd-mom}, one has to have
\[
\sum_{i=1}^\infty \alpha_i + \sum_{i=1}^\infty \beta_i \in \ensemble{0, 1}
\]
and there has to be at most two values in $\ensemble{\gamma_i \suchthat i \in \Z \setminus \ensemble{ 0} } = \ensemble{\alpha_i, \beta_i \suchthat i = 1, 2,\ldots}$, including the mandatory $0$.  On the other hand, in order to have the equality in \eqref{eq:diff-sq-smaller-sum-sq}, either $(\alpha_i)_i$ or $(\beta_i)_i$ has to be identically $0$.  This proves the implication \eqref{irred-phenom} $\Rightarrow$ \eqref{irred-par-condi}.

  \eqref{irred-par-condi} $\Rightarrow$ \eqref{irred-phenom}: We first consider the case \eqref{par-tensor-trace-case}.  In this case one can realize $M_\kappa$ as the fixed point algebra of the adjoint action of a compact group on the infinite tensor product of matrix algebras (an ITPFI action,  \citelist{\cite{MR696688}*{Section 5.3}\cite{wassermann-thesis}*{Section III.4}}) as follows.
  
  Let $V$ be a Hilbert space of dimension $n$.  The permutation of tensors
  \[
  s( v_0 \otimes \cdots \otimes v_{m-1}) = v_{s^{-1} 0} \otimes \cdots v_{s^{-1} (m-1)}
  \]
  defines a representation $\pi^{(n)}_m$ of $S_m$ on $V^{\otimes m}$.  The normalized trace $\tr$ on $B(V^{\otimes m})$ restricts to $\tau_\kappa|_{S_m}$ on $C^*S_m$ via $\pi^{(n)}_m$.  The system of embeddings
  \[
  B(V^{\otimes m}) \simeq B(V)^{\otimes m} \rightarrow B(V)^{\otimes m + 1}  \simeq B(V^{\otimes m + 1}), X \mapsto X \otimes 1
  \]
  is compatible with the embeddings $S_m \rightarrow S_{m+1}$ for $m \in \N$.
  
  There is an action $\Ad$ of the projective unitary group $\PU(V)$ of $V$ on $B(V)$ by the conjugation.  It induces the tensor product action $\Ad^{\otimes m}$ of $\PU(V)$ on $B(V)^{\otimes m}$, and the image of $C^*S_m$ is the fixed point algebra of this action.  Hence we have a compatible system
  \[
  (\Ad^{\otimes m}\colon \PU(V) \curvearrowright B(V^{\otimes m}), \pi^{(n)}_m(C^*S_m) = B(V^{\otimes m})^{\PU(V)}, \tr|_{C^*S_m} = \tau_\kappa)
  \]
  of actions of $\PU(V)$ on matrix algebras and a trace preserving homomorphism of $C^*S_m$ onto the fixed point algebras.
  
  Let $M$ denote the II$_1$ factor given as the closure of $\varinjlim B(V^{\otimes m})$ with respect to the tracial state $\tr$.  The action $\Ad^{\otimes m}$ induces an action $\Ad^{\otimes \infty}$ of $\PU(V)$ on $M$.  The above system gives an identification of $M_\kappa$ with the $\Ad^{\otimes \infty}$-fixed point subalgebra $M^{\PU(V)}$ of $M$.
  
  Moreover, $M$ admits a subfactor $N$ generated by the union of the subalgebras $1 \otimes B(V^{\otimes m - 1})$ of $B(V^{\otimes m})$ for $m \in \N_{>0}$, which is stable under the action of $\PU(V)$.  The fixed point subalgebra $N^{\PU(V)}$ of $N$, which agrees with the intersection of $M_\kappa$ and $N$, is equal to $N_\kappa$.  For such a construction one has the irreducibility $(N^{\PU(V)})' \cap M^{\PU(V)} = \C$.
  
  The case \eqref{eq:thoma-par-itpfi-perm-det} also follows from the above argument.  Indeed, the mapping $\theta\colon \sigma \mapsto (-1)^{\absolute{\sigma}} \sigma$ on $\C[S_\infty]$ defines an $*$-algebra automorphism of $C^*S_\infty$.  If the parameter $\kappa$ is of the form \eqref{eq:thoma-par-itpfi-perm-det}, let $\kappa'$ denote the parameter of the form \eqref{eq:thoma-par-itpfi-perm} obtained by exchanging $\alpha_i$ and $\beta_i$  in $\kappa$ for each $i \in \N$.  Then one has $\tau_{\kappa'} = \tau_\kappa \circ \theta$, and $\theta$ extends to an isomorphism $\sigma\colon M_{\kappa'} \rightarrow M_{\kappa}$ satisfying $N_\kappa = \theta(N_{\kappa'})$.  Hence it reduces to the case \eqref{eq:thoma-par-itpfi-perm}, and we also have $N_\kappa' \cap M_\kappa = \C$ in this case.
  
  It remains to consider the case \eqref{par-reg-rep-case}.  In this case the trace $\tau_\kappa$ is the standard trace
  \[
  \tau(\sum_{s \in S_\infty} c_s \lambda_s) = c_e
  \]
  of the discrete group $S_\infty$.  Hence the GNS representation obtained from $\tau_\kappa$ can be identified with the left regular representation of $S_\infty$.
  
  The conjugation by the elements of $S_{1\le}$ on $S_\infty$ defines an equivalence relation whose equivalence classes are infinite sets except for the trivial class $\ensemble{e}$ consisting of the neutral element.  The operators in the relative commutant of $S_{1\le}$ inside the left regular von Neumann algebra $L S_\infty$ should be represented by the functions that are constant on each equivalence class of this relation.  Hence $N_\kappa = L S_{1\le}$ is irreducible inside $M_\kappa = L S_\infty$.
\end{proof}

\begin{rmk}
 Gohm--K\"{o}stler~\cite{arXiv:1005.5726} also obtained the content of the above theorem with a different method.
\end{rmk}

\subsection{Structure of the relative commutant when \texorpdfstring{$\gamma = 0$}{gamma = 0}}\label{subsec:gamma-zero-rel-comm}

In the rest of the section we compute the relative commutant $N_\kappa' \cap M_\kappa$ for the case $\gamma = 0$ (i.e. $\sum_i \alpha_i + \sum_j \beta_j = 1$).  In such a case one has a representation of $(C^*S_\infty, \tau_\kappa)$ into the infinite tensor product $(M_N(\C), \phi_{\alpha, \beta})^{\otimes \infty}$ of matrix algebras for some fixed state $\phi_{\alpha, \beta}$ on $M_N(\C)$ ($N \in \N \cup \ensemble{\infty}$) as defined by Vershik--Kerov~\cite{MR1104274}.

Let us recall their construction.  We put
\begin{align*}
m_\alpha &= \sup  \ensemble{ i \in \N_{>0} \suchthat \alpha_i > 0}, & m_\beta &= \sup \ensemble{ j \in \N_{>0} \suchthat \beta_{j} > 0}, \quad \text{and}
\end{align*}
\begin{equation*}
X_{\alpha, \beta} = ([-m_\beta, -1] \cup [1, m_\alpha]) \cap \Z,
\end{equation*}
where we make conventions that $\sup \emptyset = 0$ and $[1, 0] = \emptyset = [0, -1]$.  Next, let $a_{\alpha, \beta}$ be the diagonal operator on $\ell^2 X_{\alpha, \beta}$ given by $a_{\alpha, \beta} \delta_i = \alpha_i \delta_i$ for $0 < i$ and $a_{\alpha, \beta} \delta_j = \beta_{-j} \delta_j$ for $j < 0$.  Then $a_{\alpha, \beta}$ is a positive trace class operator and $\phi_{\alpha, \beta}(T) = \Tr(a_{\alpha, \beta} T)$ is a state on $B(\ell^2 X_{\alpha, \beta})$, where $\Tr$ is the nonnormalized trace on $B(\ell^2 X_{\alpha, \beta})$.

Let $s \in S_n$ and $x = (x_k)_{k = 1}^n \in X_{\alpha, \beta}^n$.  We define
\[
\sign(s, x) = \# \ensemble{(k, l) \suchthat {x_k, x_l < 0, 1\le k < l \le n \text{ and } s(k) > s(l)}} \pmod 2
\]
and $\tau(s, x)\in \ensemble{1, -1}$ by
\[
\tau(s, x) = (-1)^{\sign(s, x)}.
\]
We also put $s.x$ for the permuted sequence $(s.x)_i = x_{s^{-1}(i)}$. Then $s . \delta_x = \tau(s, x) \delta_{s.x}$ defines a unitary representation of $S_n$ on $\ell^2 X_{\alpha, \beta}^{\otimes n}$.

The system $C^*S_n \rightarrow B(\ell^2 X_{\alpha, \beta})^{\otimes n}$ of homomorphisms obtained this way is compatible with the inclusions $S_n \rightarrow S_{n + 1}$ and $B(\ell^2 X_{\alpha, \beta})^{\otimes n} \rightarrow B(\ell^2 X_{\alpha, \beta})^{\otimes n + 1}$.  Hence we obtain a homomorphism $\pi_\kappa\colon C^*S_\infty \rightarrow B(\ell^2 X_{\alpha, \beta})^{\otimes \infty}$.  The pullback of the state $\phi_{\alpha, \beta}^{\otimes \infty}$ on $B(\ell^2 X_{\alpha, \beta})^{\otimes \infty}$ to $C^*S_\infty$ agrees with the trace $\tau_\kappa$.  The image of the subalgebra $C^*S_{1 \le }$ is contained in $1_{\ell^2 X_{\alpha, \beta}} \otimes B(\ell^2 X_{\alpha, \beta})^{\otimes \infty}$.

For any pair $(i, i')$ of elements in $X_{\alpha, \beta}$, let $e_{i, i'}$ denote the partial isometry $\delta_j \mapsto \delta_{i,j} \delta_{i'}$ in $B(\ell^2 X_{\alpha, \beta})$.  We also write $e_i$ to denote the projection $e_{i, i}$.

With the above constructions, one has the following refinement of Theorem~\ref{thm:irred-characterization} in this situation.

\begin{lem}\label{lem:rel-comm-projs}
 In the above notation, we have
 \[
 E^{M_\kappa}_{N_\kappa' \cap M_\kappa}((0, 1)) = \sum_{i=1}^{m_\alpha} \alpha_i e_i \otimes 1 - \sum_{j=-1}^{-m_\beta} \beta_{-j} e_{j} \otimes 1.
 \]
\end{lem}

\begin{proof}
We prove this by estimating the $2$-norm of
\begin{equation}\label{eq:diff-cond-part-condi}
\left( \sum_{i=1}^{m_\alpha} \alpha_i e_i \otimes 1 - \sum_{j=-1}^{-m_\beta} \beta_{-j} e_{j}\otimes 1\right) - E_{\pi_\kappa(T_k)' \cap M_\kappa}((0, 1))
\end{equation}
with respect to $\phi_{\alpha,  \beta}^{\otimes\infty}$ as $k \rightarrow \infty$.  By Lemma~\ref{lem:condi-exp-transp}, the square of the $2$-norm of \eqref{eq:diff-cond-part-condi} is equal to
\begin{multline}\label{eq:diff-2-norm-expand}
 \phi_{\alpha, \beta}^{\otimes k}\left(\sum_{i=1}^{m_\alpha} \alpha_i^2 e_i \otimes 1 + \sum_{j=-1}^{-m_\beta} \beta_{-j}^2 e_{j} \otimes 1\right) + \tau_\kappa\left(\frac{k-1 + \sum_{i \neq j} (0, i, j)}{(k-1)^2}\right)\\
 - \frac{2}{k-1} \sum_{l=1}^{k-1} \phi_{\alpha, \beta}^{\otimes k}\left(\left( \sum_{i=1}^{m_\alpha}\alpha_i e_i  \otimes 1 - \sum_{j=-1}^{-m_\beta}\beta_{-j} e_{j} \otimes 1\right) \pi_\kappa((0, l))\right).
\end{multline}
One has
\begin{multline*}
\phi_{\alpha, \beta}^{\otimes k}\left((e_i \otimes 1_{B(\ell^2 X_{\alpha, \beta})^{\otimes \infty}}) . \pi_\kappa((0, l))\right)\\
 = \pm \phi_{\alpha, \beta}^{\otimes\infty}\left(e_i \otimes 1_{B(\ell^2 X_{\alpha, \beta})^{ \otimes l - 1}} \otimes e_i \otimes 1_{B(\ell^2 X_{\alpha, \beta})^{\otimes \infty}}\right)
 \end{multline*}
in which the leading sign of the left hand side is determined according to that of $i$.  Hence \eqref{eq:diff-2-norm-expand} can be computed as
\[
 \frac{1}{k-1} + \left(1 + \frac{k-2}{k-1} - 2 \right) \left(\sum_{i=1}^{m} \alpha_i^3 + \sum_{j=-1}^{-n} \beta_{-j}^3\right).
\]
This converges to $0$ as $k \rightarrow \infty$.  We have thus obtained
\[
E_{N_\kappa' \cap M_\kappa}((0, 1)) = \lim_{k \rightarrow \infty} E_{\pi_\kappa(T_k)' \cap M_\kappa}((0, 1)) = \left( \sum_{i=1}^{m} \alpha_i e_i - \sum_{j=-1}^{-n} \beta_{-j} e_{j}\right),
\]
which proves the assertion.
\end{proof}

Let $P$ be the von Neumann algebra generated by $B(\ell^2 X_{\alpha, \beta})^{\otimes \infty}$ under the GNS representation with respect to the state $\phi_{\alpha, \beta}$ and let $Q$ be its sub-von Neumann algebra generated by $\C 1_{B(\ell^2 X_{\alpha, \beta})} \otimes B(\ell^2 X_{\alpha, \beta})^{\otimes \infty}$.  We then have the conditional expectation
\[
E^P_Q = \phi_{\alpha, \beta} \otimes \Id_{B(\ell^2 X_{\alpha, \beta})^{\otimes \infty}} \colon P \rightarrow Q
\]
and a square of von Neumann algebras
\begin{equation}\label{eq:dyn-sys-real-tensor-prod-rep}
\begin{CD}
Q @>>> P\\
@AAA @AAA\\
N_\kappa @>>> M_\kappa.
\end{CD}
\end{equation}

\begin{lem}\label{lem:dyn-sys-tensor-prod-rep-comm-sq}
The inclusion \eqref{eq:dyn-sys-real-tensor-prod-rep} of von Neumann algebras is a commuting square.
\end{lem}

\begin{proof}
We know that
\[
E^P_Q(\pi_\kappa((0, 1)) = 1_{B(\ell^2 X_{\alpha, \beta})} \otimes \left( \sum_i \alpha_i e_i - \sum_j \beta_{-j} e_j \right) \otimes 1_{B(\ell^2 X_{\alpha, \beta})^{\otimes \infty}}
\]
is contained in $N_\kappa$ by Lemma~\ref{lem:rel-comm-projs}, hence
\[
E^P_Q(\pi_\kappa(s_1 \cdot (0, 1) \cdot s_2)) = \pi_\kappa(s_1) E^P_Q( \pi_\kappa(0, 1)) \pi_\kappa( s_2)
\]
is in $N_\kappa$ for any $s_1, s_2 \in S_{1 \le}$.

On the other hand, we have the double coset decomposition
\[
S_\infty = S_{1 \le} \coprod (S_{1 \le} \cdot (0, 1) \cdot S_{1 \le})
\]
of $S_\infty$ with respect to $S_{1\le}$.  Hence the operator $E^P_Q(\pi_\kappa(s)) \in Q$ actually belongs to $N_\kappa$ for any $s \in S_\infty$.
\end{proof}

Consider a partition $\coprod_{k=1}^p I_k$ of the index set $[1, m_\alpha]$ such that $i$ and $i'$ lie in the same element of the partition if and only if $\alpha_i = \alpha_{i'}$.  Similarly, define a partition $\coprod_{l = -q}^{-1} J_l$ of $[-m_\beta, -1]$.
 
 For each $1 \le k \le p$ let $f_k$ denote the projection $\sum_{i \in I_k} e_i \in B(\ell^2 X_{\alpha, \beta})$.  Analogously we define $f_l = \sum_{j \in J_l} e_j$ for each $-q \le l \le -1$.

\begin{thm}\label{thm:dyn-sys-rep-min-projs-rel-comm}
 Suppose that we have $\gamma = 0$ in the Thoma parameter $\kappa$.  Then $N_\kappa' \cap M_\kappa$ is commutative, and the projections $f_k \otimes 1_{B(\ell^2 X_{\alpha, \beta})^{\otimes \infty}}$ for $-q \le k \le -1$ and $1 \le k \le p$ are the minimal projections of $N_\kappa' \cap M_\kappa$.
\end{thm}

\begin{proof}
 These projections are the spectral projections of $E_{N_\kappa' \cap M_\kappa}(\pi_\kappa((0, 1)))$ by Lemma~\ref{lem:rel-comm-projs}.  In particular, they belong to $N_\kappa' \cap M_\kappa$.  Moreover, $N_\kappa' \cap P$ is equal to $B(\ell^2 X_{\alpha, \beta}) \otimes 1_{B(\ell^2 X_{\alpha, \beta})^{\otimes \infty}}$.  Hence it remains to show that
 \[
 M_\kappa \cap (f_k B(\ell^2 X_{\alpha, \beta}) f_k) \otimes 1 = \C f_k \otimes 1
 \]
 for any integer $k$.
 
Let $n$ be an arbitrary positive integer.  As in Lemma~\ref{lem:dyn-sys-tensor-prod-rep-comm-sq}, the conditional expectation
\[
E_{B(\ell^2 X_{\alpha, \beta})} = \Id_{B(\ell^2 X_{\alpha, \beta})} \otimes \phi_{\alpha, \beta}^{\otimes n - 1}\colon B(\ell^2 X_{\alpha, \beta})^{\otimes n} \rightarrow B(\ell^2 X_{\alpha, \beta}) \otimes \C 1_{B(\ell^2 X_{\alpha, \beta})^{\otimes n-1}}
\]
with respect to the state $\phi_{\alpha, \beta}^{\otimes n}$ preserves the subalgebra $M_\kappa \cap B(\ell^2 X_{\alpha, \beta})^{\otimes n}$. Hence the map
 \[
\psi_n = B(\ell^2 X_{\alpha, \beta})^{\otimes n} \rightarrow B(\ell^2 X_{\alpha, \beta})_{f_k}, x \mapsto f_k E_{B(\ell^2 X_{\alpha, \beta})}(x) f_k
\]
restricts to the conditional expectation of $(M_\kappa)_{f_k \otimes 1}$ onto $(N_\kappa' \cap M_\kappa)_{f_k \otimes 1}$.

Then it reduces to show that, for any $n$, the image of any $s \in S_n$ under $\psi_n$ is a scalar multiple of $f_k$.  Using the matrix elements $e_{i, i'} \in B(\ell^2 X_{\alpha, \beta})$, one has
\begin{equation}\label{eq:pi-k-s-dyn-sys-expr}
\pi_\kappa(s)  = \sum_{x = (x_0, \ldots, x_{n-1}) \in X_{\alpha,\beta}^n} \tau(s, x) e_{x_0,(s.x)_0} \otimes \cdots \otimes e_{x_{n-1},(s.x)_{n-1}}.
\end{equation}

When one applies $E_{B(\ell^2 X_{\alpha, \beta})}$ to \eqref{eq:pi-k-s-dyn-sys-expr}, the only surviving terms are the ones corresponding to $x \in X_{\alpha,\beta}^n$ which are constant on the orbits of $s$.  Let $\omega_0$ denote the orbit of $0 \in \ensemble{0, \ldots, n - 1}$ under $s$ and let $\omega_1, \ldots, \omega_r$ be the other nontrivial orbits.  Suppose that $x$ is constant on each orbit $\omega_q$, and let $x_q$ be the value of $x$ on $\omega_q$.  Then one has
\[
 \tau(s, x) = \prod_{\substack{q=0,\ldots,r\\x_q < 0}}(-1)^{\#\omega_q-1}.
\]

 Thus the image of $\pi_\kappa(s)$ under $E_{B(\ell^2 X_{\alpha, \beta})}$ is equal to
\[
\sum_{(x_0,\ldots,x_r) \in X_{\alpha, \beta}^{r + 1}} \gamma_{x_0}^{\# \omega_0 - 1} \prod_{q=1}^r \gamma_{x_q}^{\# \omega_q} \prod_{\substack{q=0,\ldots,r\\x_q < 0}} (-1)^{\#\omega_q-1} e_{x_0} \otimes 1
\]
for the sequence $(\gamma_i)_{i \in \Z \setminus \ensemble{0}}$ as in the proof of Theorem~\ref{thm:irred-characterization}.  Hence we have
\[
\psi_n(s) =  (\sign(k)\gamma_{k})^{\# \omega_0 - 1} \sum_{(x_1,\ldots,x_r) \in X_{\gamma, \beta}^{r}} \: \prod_{q=1}^r (\gamma_{x_q})^{\# \omega_q} \prod_{\substack{q=1,\ldots,r\\x_q < 0}} (-1)^{\#\omega_q-1}  f_k \otimes 1,
\]
which proves the assertion.
 \end{proof}

\begin{cor}
Suppose that $1 > \alpha_{1} > \frac{1}{2}$ and $\alpha_{2} = 1 - \alpha_{1}$ in the Thoma parameter $\kappa$.  Then the inclusion $N_\kappa \subset M_\kappa$ is isomorphic to the inclusion of von Neumann algebras $N \subset M$ generated by Jones projections $M = \langle e_0 , e_1, \ldots, \rangle'', N = \langle e_1 , e_2, \ldots, \rangle''$ satisfying $e_i e_{i\pm 1} e_i = \delta e_i$ for $\delta = \alpha_{1} \alpha_{2}$ and $[e_i, e_j] = 0$ for $\absolute{i-j} > 1$.
\end{cor}

\begin{proof}
As in the above description one finds a state $\psi_0$ on $M_2^{\otimes \infty}$ and a compatible system of $*$-homomorphisms $\phi_n\colon C^*S_n \rightarrow M_2^{\otimes n}$ which together induce a $*$-homomorphism $\phi\colon C^*S_\infty \rightarrow M_2^{\otimes \infty}$, such that $\psi = \psi_0^{\otimes \infty}$ on $M_2^{\otimes \infty}$ restricts to $\tau_\kappa$.  The state $\psi_0$ can be written as $x \mapsto \Tr(x a)$, where $\Tr$ is the nonnormalized trace and $a$ is a diagonal matrix with positive entries $\alpha_1$ and $\alpha_2$.  Notice that this is a particular case of the well known constructions as found in \cite{MR680655}*{Theorem 4.2} and \cite{wassermann-thesis}*{Theorem III.4}.
 
 The centralizer algebra $M = (M_2^{\otimes \infty})''^{\psi}$ of $\psi$ and its subalgebra $N = (M_2 \otimes 1_{M_2^{\otimes \infty}})' \cap M$ is indeed isomorphic to the subfactor inclusion coming from the Jones projections (Pimsner and Popa \cite{MR860811}*{Section 5.5}).  Now, the image of $C^*S_\infty$ is clearly contained in $M$, while that of $C^*S_{1\le}$ is in $N$.  In fact, the image of $C^*S_\infty$ is equal to $M$ and so is that of $C^*S_{1\le}$ to $N$: since $M_2 \cap M_\kappa = N_\kappa' \cap M_\kappa$ is not trivial, it contains the minimal projections $e_{11}$ and $e_{22}$ of $M_2$ corresponding to the diagonalization of $\psi_0$.  If we denote the corresponding minimal projections in the $n$-th component in the tensor product $M_2^{\otimes \infty}$ by $e^n_{ii}$ and the corresponding partial isometries by $e^n_{ij}$, the transposition $(n, n+1)$ is represented by 
\[
x_n = e^n_{11} e^{n+1}_{22} + e^n_{22} e^{n+1}_{11} + e^n_{12} e^{n+1}_{21} + e^n_{21} e^{n+1}_{12}.
\]
 Hence $M_\kappa$ contains $e^n_{11} x_n e^n_{22} = e^n_{12} e^{n+1}_{21}$, etc., which allows us to realize the representative
 \[
 \alpha_{2} e^n_{11} e^{n+1}_{22} + \alpha_{1} e^n_{22} e^{n+1}_{11} + \sqrt{\alpha_{1} \alpha_{2}} (e^n_{12} e^{n+1}_{21} + e^n_{21} e^{n+1}_{12})
 \]
 inside $M_\kappa$ of the Jones projections corresponding to the index $(\alpha_{1} \alpha_{2})^{-1}$.
 \end{proof}

\begin{rmk}\label{rmk:af-ladder}
The group algebra $C^*S_\infty$ is an AF-algebra, being the inductive limit of the increasing sequence of the algebras $C^*S_n$.  The subalgebra $C^*S_{1\le}$ is also the inductive limit of the subalgebras $C^*T_n$ of $C^*S_n$.  Hence we obtain squares of finite dimensional algebras
 \begin{equation}\label{eq:subf-bratelli-diag}
 \begin{CD}
 \cdots @>>> C^*T_n @>>> C^*T_{n+1} @>>> \cdots \\
  @. @VVV @VVV @. \\
 \cdots @>>> C^*S_n @>>> C^*S_{n+1} @>>> \cdots .
 \end{CD}
 \end{equation}
 
 The cases where these squares become commuting (i.e. $[E^{C^*S_{n+1}}_{C^*S_n}, E^{C^*S_{n+1}}_{C^*T_{n+1}}] = 0$) with respect to the trace $\tau_\kappa$ happens to agree with the ones of Theorem~\ref{thm:irred-characterization}.  Indeed, if the above square is commuting for $n = 2$, the image of the transposition $(1\ 2) \in C^*T_3$ under $E^{C^*S_3}_{C^*S_2}$ should lie in $\C = C^*T_2$.  If that is the case, we will have $\tau_\kappa((0\ 1)(1\ 2)) = \tau_\kappa(0\ 1)\tau_\kappa(1\ 2)$, which is equivalent to the formula \eqref{eq:s-tau-square-eq-3-elem-tau}.
 \end{rmk}

\section{Index of \texorpdfstring{$S_\infty$}{S-infinity} subfactors}

\begin{lem}\label{lem:ep-and-k-fac-diminish}
 Let $R > 0$ be an arbitrary positive real number.   The sequence of numbers
 \[
 \frac{R^{k(k+1)/2} k! (k-1)! \cdots 1!}{(k(k+1)/2)!}\quad (k \in \N)
 \]
 converges to $0$ as $k \rightarrow \infty$.
\end{lem}

\begin{proof}
The natural logarithm of the above expression is equal to
\begin{equation}\label{eq:log-dim-irred-w-eps}
 \frac{k(k+1)}{2} \log R + \sum_{i=1}^k \sum_{j=1}^i \log j - \sum_{i=1}^{k(k+1)/2} \log i.
\end{equation}
We can estimate 
\[
 \sum_{j=1}^a \log j \ge \int_1^{a} \log(x) d x = a\log(a) - a
\]
 and
\[
 \sum_{i=1}^a \sum_{j=1}^i \log j \le \int_2^{a+1} \left. (x+1)(\log(x+1) - 1) \right. d x \le \frac{(a+2)^2 \log(a+2)}{2} - \frac{3 (a+2)^2}{4}
\]
for any positive integer $a$.  Hence \eqref{eq:log-dim-irred-w-eps} has the estimate from above by
\begin{multline*}
 \frac{k(k+1)}{2} \log R + \frac{(k+2)^2}{2} \log(k+2) - \frac{3 (k+2)^2}{4} \\
 - \frac{k(k+1)}{2} \log\left(\frac{k(k+1)}{2}\right) + \frac{k(k+1)}{2}.
\end{multline*}
This can be reorganized as the sum of
\begin{gather*}
\frac{k(k+1)}{2} \left( (\log(R)+1) -  \frac{\log(k)}{2}\right),\\
\frac{3k+4}{2} \log(k+2) - \frac{3(k+2)^2}{4},
\intertext{and}
\frac{k(k+1)}{2} \left( \log(k+2) - \log\left(\frac{\sqrt{k}(k+1)}{2}\right)\right).
\end{gather*}
All of the above go to $-\infty$ as $k \rightarrow \infty$, which shows that \eqref{eq:log-dim-irred-w-eps} converges to $-\infty$.
\end{proof}

\begin{lem}\label{lem:small-proj}
 Let $\tau$ be a faithful tracial state on $C^*S_\infty$.  For any positive real number $\epsilon$, there exists a positive integer $n$ and a projection $e \in C^*S_n$ such that $0 < \tau(e) < \epsilon^n$.
\end{lem}

\begin{proof}
 We consider the Young diagrams with the isosceles right triangle shape, i.e. the ones having rows of length $k, k-1, \ldots, 1$ exactly once for each for $k \in \N$, with the total number of boxes equal to $k (k + 1) / 2$.  By the hook-length formula, the irreducible representation of $S_{k (k + 1) / 2}$ corresponding to this diagram has dimension
 \[
 \frac{(k (k + 1) / 2)!}{3^{k-1}5^{k-2}\cdots (2k-1)}.
 \]
 Hence a minimal projection $e_k$ belonging to the factor of this representation satisfies
 \[
0 < \tau(e_k) < \frac{3^{k-1}5^{k-2}\cdots (2k-1)}{(k (k + 1) / 2)!} < 2^{k (k + 1) / 2} \frac{k! (k-1)! \cdots 1!}{(k(k+1)/2)!}.
 \]
 Applying Lemma~\ref{lem:ep-and-k-fac-diminish} to $R = 2 \epsilon^{-1}$, we obtain
\[
 \left(2\epsilon^{-1}\right)^{k (k + 1) / 2} \frac{k! (k-1)! \cdots 1!}{(k(k+1)/2)!} < 1
\]
 for large enough $k$.  For such $k$, taking $n = k(k+1)/2$ and the corresponding projection $e = e_k$ in $C^*S_n$ as above, we have $\tau(e) < \epsilon^n$.
\end{proof}

\begin{thm}\label{thm:subf-infin-ind}
The subfactor $N_\kappa \subset M_\kappa$ has infinite index if and only if $\tau_\kappa$ is faithful on $C^*S_\infty$.
\end{thm}

\begin{rmk}\label{rmk:fathful-trace-thoma-par}
 The trace $\tau_\kappa$ is faithful in the following cases: 1) $\alpha_i > 0$ for any $i$, 2) $\beta_i > 0$ for any $i$, and 3) $\gamma > 0$.  Otherwise it is not, as the corresponding infinite Young tableaux can be taken to have a bounded number of rows and columns~\cite{wassermann-thesis}*{Theorem III.5}.
\end{rmk}

\begin{proof}[Proof of the theorem]
Suppose that the trace $\tau_\kappa$ is faithful.  For each $i \in \N$, let $ L_{\tau_\kappa}S_{i\le}$ be the subalgebra of $M_\kappa$ generated by $S_{i \le}$.  Since there is an isomorphism $\sigma^n\colon S_\infty \rightarrow S_{n\le}$ satisfying $\sigma^n(S_{1\le}) = S_{n+1\le}$ and $\tau_\kappa(s) = \tau_\kappa(\sigma^n(s))$, the inclusions $L_{\tau_\kappa}S_{i+1\le} \subset  L_{\tau_\kappa}S_{i\le}$ are all isomorphic to $N_\kappa \subset M_\kappa$.  In particular, we have $[M_\kappa : L_{\tau_\kappa} S_{i \le}] = [M_\kappa : N_\kappa]^i$.  We also note that the image of $C^*S_n$ in $M_\kappa$ is contained in the relative commutant of $L_{\tau_\kappa} S_{n\le}$ for any $n \in \N_{>0}$.

Let $\epsilon$ be an arbitrary positive real number.  By Lemma~\ref{lem:small-proj}, there exists an integer $n$ and a projection $e \in C^*S_n$ satisfying $\tau_\kappa(e) < \epsilon^n$.  This means that the image $E_{L_{\tau_\kappa}S_{n \le}}(e)$ of $e$ under the conditional expectation onto $L_{\tau_\kappa} S_{n \le}$ is a positive scalar smaller than $\epsilon^n$.   By the Pimsner--Popa inequality, we have $[M_\kappa : L_{\tau_\kappa} S_{n \le}] > \epsilon^{-n}$.  Thus $[M_\kappa : N_\kappa] > \epsilon^{-1}$ for any $\epsilon$, which shows $[M_\kappa : N_\kappa] = \infty$ when $\tau_\kappa$ is faithful.

Next, suppose that $\tau_\kappa$ is not faithful.  By Remark~\ref{rmk:fathful-trace-thoma-par}, this implies that $\gamma = 0$ and only finitely many terms in the sequences $(\alpha_i)_{i \in \N_{>0}}$ and that $(\beta_j)_{j\in\N_{>0}}$ are nonzero.  Hence we can apply the construction of Section~\ref{subsec:gamma-zero-rel-comm} to this situation, and there are algebras $Q \subset P$ and a commuting square as in \eqref{eq:dyn-sys-real-tensor-prod-rep} by Lemma~\ref{lem:dyn-sys-tensor-prod-rep-comm-sq}.  Again by the Pimsner--Popa inequality, in order to conclude $[M_\kappa : N_\kappa] < \infty$ it is enough to show that there is a constant $C > 0$ satisfying $E^P_Q(y) \ge C y$ for any $y \in P_+$.  

Put $\delta = \min (\alpha_{m_\alpha}, \beta_{m_\beta})$ so that we have $\phi_{\alpha, \beta} - \delta \tr \ge 0$ on $B(\ell^2 X_{\alpha, \beta})$. Then we have the estimate
\[
E^P_Q(y) \ge \delta \tau \otimes \Id (y) \ge \delta (m_\alpha + m_\beta)^{-2} y
\]
for $y \in P_+$.  This shows that $[M_\kappa : N_\kappa] < \infty$ when $\tau_\kappa$ is not faithful.
\end{proof}

\subsection{Relative Entropy}

In the rest of the paper we consider the relative entropy of the inclusion $N_\kappa \subset M_\kappa$.  Recall that it is defined by
\[
H_\tau(M_\kappa \mid N_\kappa) = \sup_{\substack{k \in \N, x_j \in (M_{\kappa})_+,\\ x_1+ \cdots + x_k = 1}} \sum_{j=1}^k \tau_\kappa(\eta(E_{N_\kappa}^{M_\kappa}(x_j))),
\]
where $\eta(t) = -t \log(t)$.  We follow the conventions of Neshveyev--St{\o}rmer \cite{MR2251116}*{Chapter 10} in the following.

\begin{thm}\label{thm:subf-entropy-estim}
 We have the estimate of the relative entropy 
\begin{equation}\label{eq:rel-entropy-l-le-r}
H_{\tau_\kappa}(M_\kappa \mid N_\kappa) \le 2 \sum_{i} \eta(\alpha_i) + \eta(\beta_i). 
\end{equation}
The equality holds when $\gamma = 0$ and there is no duplicate in each of the sequences $(\alpha_i)_{i \in \N}$ and $(\beta_j)_{j \in \N}$ other than $0$.
\end{thm}

\begin{proof}
 We are going to use the AF-structure of the inclusion $\sigma(M_\kappa) = N_\kappa \subset M_\kappa$ given by the squares \eqref{eq:subf-bratelli-diag}.   Although they are not commuting in general, we still have 
  \[
  H_{\tau_\kappa}(M_\kappa \mid N_\kappa) \le \varliminf H_{\tau_\kappa}(C^*S_n \mid C^*T_n)
  \]
 by \cite{MR860811}*{Proposition 3.4}.
 
Let $G$ be the group of finitely supported permutations of the set $\Z$.  Its group algebra admits `the shift automorphism' $\sigma$ characterized by $\sigma\colon (n, m) \rightarrow (n+1, m+1)$ on the transpositions.  For any Thoma parameter $\kappa$, the trace $\tau_\kappa$ on $C^*S_\infty$ admits a unique extension as a trace (which we still denote by $\tau_\kappa$) to $C^*G$.

 The subalgebras $\sigma^{kn}(C^*S_n)$ of $C^*G$ for $k \in \N$ are $\tau_\kappa$-independent by the multiplicativity of $\tau_\kappa$, e.g. $\tau_\kappa(x \sigma^n(y)) = \tau_\kappa(x) \tau_\kappa(y)$ for $x, y \in C^*S_n$.  By \cite{MR2251116}*{Proposition~10.4.5},  the sequence of subalgebras $(C^*S_n)_{n \in \N}$ is a generating sequence for $\sigma$.  In particular, we can compute the dynamical entropy of $\sigma$ as
  \[
  h_{\tau_\kappa}(\sigma) = \lim_{n\rightarrow\infty}\frac{1}{n}H_{\tau_\kappa}(C^*S_n).
  \]
  
  Now, \cite{MR2251116}*{Lemmas 10.3.3 and 10.3.4} apply to our case without change.  We then have
  \begin{equation*}%\label{eq:cesaro-entropy-shift-center-entropy}
  \frac{1}{N}\sum_{n=1}^N H_{\tau_\kappa}(C^*S_n \mid C^*T_n) = \frac{2}{N}H_{\tau_\kappa}(C^*S_N) -  \frac{1}{N}H_{\tau_\kappa}(Z(C^*S_N)) + \frac{1}{N}\sum_{n=1}^NC_n
  \end{equation*}
  with positive numbers $C_n$ satisfying
  \[
  \frac{1}{N}\sum_{n=1}^NC_n \rightarrow 0 \quad (N \rightarrow \infty).
  \]
  
  Hence it follows that $\varliminf H_{\tau_\kappa}(C^*S_n \mid C^*T_n) \le 2 h_{\tau_\kappa}(\sigma)$.  On the other hand, by the result of Boyko--Nessonov~\cite{MR2172640}, we have
 \[
 h_{\tau_\kappa}(\sigma) = \sum \eta(\alpha_i) + \sum \eta(\beta_i)
 \]
with respect to $\tau_\kappa$, which proves \eqref{eq:rel-entropy-l-le-r}.

It remains to show the converse inequality to \eqref{eq:rel-entropy-l-le-r} when the conditions $\gamma = 0$, $\alpha_i = \alpha_{i'} \Rightarrow i = i'$ or $\alpha_i = 0$, and $\beta_j = \beta_{j'} \Rightarrow j = j'$ or $\beta_j = 0$ hold.

By Theorem~\ref{thm:dyn-sys-rep-min-projs-rel-comm} and the assumption, the projections $e_{i,i}$ for $-q \le i \le -1$ and $1 \le i \le p$ are the minimal projections of $M_\kappa \cap N_\kappa'$.  Then \cite{MR860811}*{Theorem 4.4} implies that
\[
H_{\tau_\kappa}(M_\kappa | N_\kappa) = 2 \sum \eta(\tau_\kappa(f_k)) + \sum \tau_\kappa(f_k) \log([(M_\kappa)_{f_k}: (N_\kappa)_{f_k}]).
\]
Combined with $\tau_\kappa(f_k) = \alpha_k$ when $k > 0$ and $\tau_\kappa(f_k) = \beta_{-k}$ when $k < 0$, we obtain the desired equality.
\end{proof}

\paragraph{Acknowledgements} The author would like to thank R. Tomatsu for his suggestion to work on this subject and for many encouraging discussions with him in the early stages of the project.  He is also grateful to S. Neshveyev, T. Katsura, Y. Oshima, and the referee for many fruitful comments on various occasions.  The author also benefited from conversations with A. Zuk, Y. Ueda, E. Blanchard, A. Hora, D. Ara, M. Hajli, M.~Imsatfia, S. Oguni, T. Fukaya, T. Kato, M. Mimura, and N. Ozawa.  Last but not least, he is deeply indebted to Y. Kawahigashi and G. Skandalis for their continuous support throughout the period of this research.

% \bibliographystyle{math}

\end{document}